\documentclass[letterpaper]{physor2026}

\usepackage{caption}

\usepackage{algorithm}
\usepackage{algpseudocode}
\usepackage{afterpage}
\usepackage{amsbsy}
\usepackage{amsmath}
\usepackage{amssymb}
\usepackage{amsthm}
\usepackage{bm}
\usepackage{bbm}
\usepackage{cases}
\usepackage{color}
\usepackage{dcolumn}
\usepackage{dsfont }
\usepackage{fonts}

\usepackage{microtype}
\usepackage{multirow}
\usepackage[margin=0.75cm,belowskip=5pt,aboveskip=1pt]{subcaption}
\usepackage{tabularx}
\usepackage{tikz}
\usepackage{url}
\usepackage{xparse}
\usepackage{physics}\usepackage{nomencl} 

\makenomenclature

\title{Scatter-Limited Hybrid Monte Carlo, Deterministic Transport with Quasi–Monte Carlo Sampling}

\addAuthor[jkrotz@nd.edu]{Johannes Krotz}{1}
\addAuthor{Ryan G. McClarren}{1}

\addAffiliation{1}{Department of Aerospace and Mechanical Engineering, University of Notre Dame, Notre
				Dame, IN}

\newcommand{\x}{{\bm{x}}}

\newcommand{\del}{\partial}
\newcommand{\Om}{{\bm{\Omega}}}
\newcommand{\sigmat}{\sigma_\mathrm{t}}
\newcommand{\sigmaa}{\sigma_\mathrm{a}}
\newcommand{\sigmas}{\sigma_\mathrm{s}}

\newcommand{\lambdas}{\lambda_\mathrm{s}}

\newcommand{\Psisub}[1]{\Psi_{\lceil #1 \rceil}}
\newcommand{\Psisup}[1]{\Psi_{\lfloor #1 \rfloor}}

\Abstract{%
We present a hybrid method for time-dependent particle transport that combines Monte Carlo (MC) estimation with a deterministic discrete ordinates (\(S_N\)) solve, augmented by quasi-Monte Carlo (QMC) sampling. For spatial discretizations, the MC component computes a piecewise-constant (cell-averaged) solution, while the \(S_N\) stage employs bilinear discontinuous finite elements. By hybridizing the formulation, the MC subproblem after a prescribed scatter limit becomes scattering-free, yielding a simple and efficient streaming/attenuation procedure. Between time steps, a simple scatter-free MC step is run to relabel the $S_N$ solution as an MC solution.

A key feature of the approach is a tunable parameter \(N_{s}\) that controls how many material collisions are handled in the (Q)MC leg before handing off to the deterministic \(S_N\) solve; \(N_s=0\) recovers a purely uncollided MC leg, while \(N_s>0\) produces multi-scatter hybrids. QMC replaces pseudorandom draws with low-discrepancy points in the existing MC sampling maps, enabling a plug-in adoption within the standard MC code with modest, localized changes.

We observe significant accuracy and convergence rate improvements through the use of QMC and practically no additional computational cost, which are generally not seen in comparable non-hybrid solves. We believe the multi-scatter approach provides additional flexibility in terms of parallelization and the choice of deterministic solver.

}


\shortTitle{PHYSOR 2026 Scatter-Limited Hybrid MC-SN with QMC Sampling}

\authorHead{J.K. \& R.G.M.}

\begin{document}

\section{Introduction}\label{sec:1}
Neutron transport is a crucial phenomenon in reactor physics. Deterministic discretizations of the full phase space are accurate but costly in seven dimensions (3 space, 2 direction of flight, 1 Energy, 1 time) \cite{lewis1984computational}, whereas Monte Carlo (MC) handles complex geometries and continuous physics but carries \(N_p^{-1/2}\) sampling error \cite{spanier2008monte}. Hybrid strategies seek the best of both: regional--coupling and biasing frameworks \cite{mosher2006incident,wagner2014fw}, and high--order/low--order (HOLO) closures \cite{park2014moment}. A complementary line exploits \emph{streaming--scattering splitting}, treating potentially anisotropic streaming separately from the more isotropic scattered field \cite{krotz_2024,hauck2013collision,crockatt2019hybrid}. Related discontinuous Galerkin (DG)--DG hybrids have demonstrated multigroup coarsening for collided components \cite{whewell2023multigroup}.

\paragraph{Contributions.}
This paper introduces a time-dependent MC--DG hybrid that: (i) advances particles with MC \emph{up to a fixed cap of \(N_s\) material scatterings} per step (the MC leg is not scatter-free for $N_s >0$); (ii) transfers particles exceeding \(N_s\) to a collided equation solved by discrete ordinates with a \emph{discontinuous Galerkin (DG)} spatial discretization; and (iii) applies a stepwise remap that resamples from the DG solution back into MC to prevent exponential depletion of low-scatter particles. With an appropriate DG scheme, the method retains the diffusion limit \cite{larsen1987asymptotic}. Numerical experiments show that this hybrid approach sees improved convergence rates by a simple replacement of pseudo-random numbers with quasi-random numbers. The bounded multi-scatter MC mechanism is new and developed in the next section. We focus on single-group MC--DG.

\section{Basics of the QMC--Hybrid Method \label{sec:2}}
	\subsection{Transport equation}
	Let $X \subset {\mathbb{R}^3}$ be a spatial domain with Lipschitz boundary and let ${\mathbb{S}^2}$ be the unit sphere in ${\mathbb{R}^3}$.  Let $\Psi= \Psi(\x,\Om,t)$ be the angular flux depending on the position $\x=(x,y,z)\in X$, the direction of flight $\Om \in \mathbb{S}^2$ and time $t>0.$ We assume that $\Psi$ is governed by the linear transport equation
	\begin{equation}\label{eq:boltzmann}
	\frac{1}{c} \del_t \Psi+\Omega\cdot \nabla_{\x}\Psi +\sigmaa \Psi = \sigmas\mathcal{S}[\Psi]   +Q,  \qquad \x \in X , \quad \Om \in {\mathbb{S}^2}, \quad  t >0,
    \end{equation}
    where  $\sigmas = \sigmas(\x)$ and $\sigmaa = \sigmaa(\x)$ are the scattering and absorption cross-sections of the material, respectively; $Q=Q(\x,\Om,t)$ is a known particle source; and $\mathcal{S}$ denotes the scattering kernel, which throughout this paper will be:
    \begin{equation}
    \mathcal{S}[\Psi] = \frac{\langle \Psi \rangle}{4\pi} -\Psi = \frac{1}{4\pi}\int\limits_{\mathbb{S}^2}\Psi \,d\Om - \Psi.  
    \end{equation}
    Note that isotropic scattering is not a restriction of the method. Indeed, this hybrid approach could be readily extended to that scattering in future work.
    
    The constant $c >0$ is the particle speed; we assume that $c=1$  for the remainder of this work.  
    The transport equation \eqref{eq:boltzmann} is equipped with initial conditions 
   \begin{equation}
   \label{eq: IC}
	\Psi(\x,\Om,0) = \tilde{\Psi}(\x,\Om), \quad\x \in X , \quad \Om \in {\mathbb{S}^2},
	\end{equation}
    and boundary condition
    \begin{equation}
    \label{eq: BC}
	\Psi(\x,\Om,t) =  G(\x,\Om,t) \qquad   \x \in \partial X, \quad \Om \cdot \bm{n}(\x) <0 ,
	\end{equation}
    where $\Psi_0$ and $G$ are known and $\bm{n}(\x)$ is the unit outward normal at $\x \in \partial X$.

	\subsection{The hybrid method}
	The hybrid method is based on a collision source splitting \cite{alcouffe2006first}. 
 Let $\Psi=\sum_{n=0}^\infty \Psi_n$, where the \textit{n-collision flux} $\left\{\Psi_n\right\}_{n=0}^\infty$ satisfy the following system of equations 
  
      \addtocounter{equation}{+1}
	\begin{align}   
	\del_t \Psi_n+\Omega\cdot \nabla_{\x}\Psi_n +\sigmaa \Psi_n &=  Q_n - \sigmas \Psi_n,  \qquad n\in \mathbb{N}_0, \label{eq:n-split}\tag{\theequation.n}   
	\end{align}
    with source terms $Q_0 = Q$ and $Q_n =\frac{\sigmas}{4\pi}\langle \Psi_{n-1} \rangle$ for all $n>0$. For simplicity we assume that $\Psi_0(\x,\Om,0)=\tilde{\Psi}(\x,\Om)$ and $\Psi_0(\x,\Om,t)=G(\x,\Om,t)$ on $\x\in\del X$ for incoming directions. Thus the initial and boundary data for $n>0$ end up being 0. This split divides particles into sets with densities $\Psi_n$, in which each particle has currently undergone $n$ collisions since its birth at the source $Q$ or entry into the domain. This further splits the scattering operator $\sigmas S[\Psi_n]$ into an absorption term $-\sigmas\Psi_n$ and a reemission term $\frac{\langle \sigmas\Psi_n \rangle}{4 \pi}$. The absorption parts shows up as a sink in the equation for $\Psi_n$, which is interpreted as a particle with $n$ collisions disappearing as particle with $n$ collisions, when it scatters an additional time, effectively being absorbed. While the isotropic reemission appears in the equation for $\Psi_{n+1}$ indicating that an $n$-scattered particle is remitted with $n+1$ collision after an additional scattering event.

Due to the linearity of \eqref{eq:boltzmann}, the splitting in \eqref{eq:n-split} is exact; that is, if all $\Psi_n$  solve their respective \eqref{eq:n-split}, then $\sum_{n=0}^\infty \Psi_n$ solves \eqref{eq:boltzmann}. A hybrid method will solve equations \eqref{eq:n-split} in two batches by splitting 
$\Psi = \sum_{n=0}^\infty \Psi_n = \Psi_0 + \sum_{n=1}^\infty \Psi_n =: \Psi_u + \Psi_c$, into the densities of uncollided and collided particles. We build on this idea by splitting this some along a different $N_s\neq 0$ into

$\Psi = \sum_{n=0}^{N_s-1} \Psi_n + \Psi_{N_s} + \sum_{n={N_s}+1}^\infty \Psi_n =: \Psi_{\lceil{{N_s}}\rceil} +\Psi_{N_s} + \Psi_{\lfloor {N_s} +1 \rfloor}$ with the \textit{sub-${N_s}$ collision flux} $\Psi_{\lceil{{N_s}}\rceil}$, the \textit{${N_s}$-collision flux}  and the \textit{super-${N_s}$-collision flux}  $\Psi_{\lfloor {N_s} +1 \rfloor}$  solving 

    \begin{subequations}
     \label{eq:split}
	\begin{align}
	\del_t \left(\Psi_{\lceil {N_s} \rceil}+\Psi_{N_s}\right)+\Omega\cdot \nabla_{\x}\left(\Psi_{\lceil {N_s} \rceil}+\Psi_{N_s}\right) +\sigmaa \Psisub{{N_s}}+\sigmat \Psi_{{N_s}} &=Q +\sigmas \mathcal{S} \left[\Psi_{\lceil {N_s} \rceil}\right]  \label{eq:uncollided} 
    \tag{\theequation.$\lceil {N_s} \rceil$}\\
	\del_t \Psi_{\lfloor {N_s}+1 \rfloor}+\Omega\cdot \nabla_{\x}\Psi_{\lfloor {N_s}+1 \rfloor} +\sigmat \Psi_{\lfloor {N_s}+1 \rfloor} &=   \sigmas\mathcal{S}\left[\Psi_{\lfloor {N_s} \rfloor} \right]  + \frac{\sigmas}{4\pi}\langle \Psi_{{N_s}}\rangle \tag{\theequation.$\lfloor {N_s}+1 \rfloor$}
 \label{eq:collided}
	\end{align}
 \end{subequations}
 
 This is a split into particles that scattered at most ${N_s}$ times and particles that scattered at least ${N_s}+1$ times. For ${N_s}=0$ the sum in $\Psisub{{N_s}}$ is empty and previous hybrid methods are recovered.
 Note that a solution to \eqref{eq:uncollided} consists of the pair $(\Psisub{{N_s}},\Psi_{N_s})$. Further note that for ${N_s}>0$ this, on its own, is not a well-posed problem and allows for infinitely many solutions. Throughout this manuscript this will be resolved by conditioning $\{\Psi_n\}_{0}^{N_s}$ to solve \eqref{eq:n-split}. This determines all $\Psi_n$ and thus $\Psisub{N_s}$ $\Psi_{N_s}$ uniquely. A short discussion of this is included in Appendix \ref{app:uniqueness}.

 In practice, \eqref{eq:uncollided} and \eqref{eq:collided} are solved at different resolutions or even with different methods.  
 Typically \eqref{eq:uncollided} is solved with a method that has high resolution in angle, and because \eqref{eq:uncollided} has no (${N_s}=0$), or lower (${N_s}>0$ but finite) coupling  in angle, it is easier to solve than \eqref{eq:boltzmann} and also easier to solve in parallel.  Meanwhile \eqref{eq:collided} inherits the full angular coupling in \eqref{eq:boltzmann}, but typically requires less angular resolution.  Within this work we will assume that \eqref{eq:uncollided} is solved using a (Q)MC--simulation, described in the next section ,and \eqref{eq:collided} is solved using a an $S_N$--solver, for example described in \cite{krotz_2024}.

Previous hybrid solvers for this system have, to our knowledge, only covered the case of ${N_s}=0$ corresponding to a split into uncollided particles and particles having undergone at least 1 scattering event. In this special case $\Psisub{{N_s}}\equiv 0$, simplifying \eqref{eq:uncollided}.
In the version discussed going forward, we will present a hybrid that can solve these equations for arbitrary ${N_s}$. This allows for a tunable division of workloads between the different methods used for each part.

    Assuming \eqref{eq:uncollided} is solved at high accuracy, while \eqref{eq:collided} is done cheaply, a majority of the workload should be carried out in the former of the two to retain high accuracy overall. However, in the non-trivial case with $\sigmas \neq 0$ particles scatter at an exponential rate, thus particle densities are, over time, transferred out of $\Psi_n$ into $\Psi_{n+1}$ for every $n$. Thus particle densities will be transferred into the $N_s$-tail flux at an exponential rate,  thus making the accuracy at which \eqref{eq:collided} is solved the driving factor in the overall accuracy.  We briefly give bounds on this transfer rate in Appendix \ref{app:transfer}

 This effect can be mitigated by abusing the autonomous nature of the equations and periodically relabelling the high scatter densities as uncollided at every time step.
 This leads to the following algorithm to propagate a solution $\Psi[t_n](\x,\Om) = \Psi(\x,\Om,t_n)$ along the time step $[t_n,t_{n+1}]$ to $\Psi[t_{n+1}]$:
 \begin{enumerate}
     \item With initial data $\Psisub{N_s}^{MC}[t_n]=\Psi[t_n]$ and $\Psi^{MC}_{N_s}\equiv 0$ solve \eqref{eq:uncollided} for solutions $\Psisub{N_s}^{MC}[t_{n+1}]$ and $\Psi_{N_s}^{MC}[t_{n+1}]$ via the scatter-limited MC algorithm with $N_s$ scatters described in the next section. Apply source $Q[t]$ and boundary data $G[t]$ in this step.
     \item Solve \eqref{eq:collided} for an $S_N$--solution $\Psi^{SN}[t_{n+1}]$ with zero initial and boundary data.
     \item (Remap) With zero allowable scatters in the MC--algorithm solve 
     \begin{align}
           \tag{\theequation.R}
\del_t \Psi_R+\Omega\cdot \nabla_{\x}\Psi_R +\sigmat \Psi_R &= \left(
Q +\sigmas \langle\Psi_{\lceil {N_s} \rceil}\rangle^{MC} \right) +\sigmas\langle \Psi_{\lfloor {N_s} \rfloor} \rangle^{MC}  + \frac{\sigmas}{4\pi}\langle \Psi_{{N_s}}\rangle^{SN} \label{eq:relabel}
     \end{align}
     for $\Psi_R^{MC}[t_{n+1}]$. Use initial data $\Psi_R[t_n] = \Psi[t_n]$ and apply boundary data $G[t]$.
  \item Set $\Psi[t_{n+1}] = \Psi_R^{MC}[t_{n+1}]$.   
 \end{enumerate}
Note that $\Psi^{MC}$ will generally be represented by a collection of computer particles distributed proportional to $\Psi^{MC}$, $\langle \Psi \rangle^{MC/SN}$ are MC and $S_N$ estimators of $\langle \Psi \rangle$ respectively, in our case an implicit capture track length estimator in the MC case and Gauss-Legendre integration of the ordinates for $S_N$.

Further note that previous hybrids replaced steps 3 and 4 with 
\begin{enumerate}
    \item[3'.]   With zero allowable scatters in the MC--algorithm solve 
     \begin{align}
           \tag{\theequation.R'}
\del_t \Psi_R+\Omega\cdot \nabla_{\x}\Psi_R +\sigmat \Psi_R &= \sigmas\langle \Psi_{\lfloor {N_s} \rfloor} \rangle^{MC}  + \frac{\sigmas}{4\pi}\langle \Psi_{{N_s}}\rangle^{SN} \label{eq:relabelprime}
     \end{align}
     for $\Psi_R^{MC}[t_{n+1}]$. Boundary and initial data in this case are zero.
    \item[4'.] Set $\Psi[t_{n+1}] =\Psisub{N_s}^{MC}[t_{n+1}]+\Psi_{N_s}^{MC}[t_{n+1}]+ \Psi_R^{MC}[t_{n+1}]$.
\end{enumerate}
These approaches are, in theory, equivalent due to the linearity of the equations involved. However, we found that the version with 3' and 4' does not benefit from QMC. Further, step 3' requires the addition of new computer particles, whereas 3 can be solved by resetting the particles from 1, thus allowing for much easier control of overall particle numbers. In practice 3 can be run by resetting particles that were added through the source terms, while particles from initial data and boundary terms can just be copied from 1 without additional computation.

\subsection{Scatter-limited Monte Carlo (MC) Algorithm}
\label{sec:scatter_limited_mc}
The algorithm described in the previous section requires an (Q)MC solver for \eqref{eq:uncollided} that (a) will scatter at most $N_s$ times and (b) tracks $\langle \Psisub{N_s} \rangle$ and $\langle\Psi_{N_s}\rangle$ separately. First we will briefly discuss the unlimited scatters case ($N_s = \infty$ and $\Psi_{N_s} = \Psi_\infty \equiv 0$), as well as the zero scatter case ($N_s = 0$ and $\Psisub{N_s}=\Psisub{0}\equiv 0$)
In all cases we use implicit capture. Thus we start by sampling particles with birth location and time proportional to \(q(x,\Omega,t)=Q(x,\Omega,t)\big/\!\int_{X\times S^2\times I}Q\,\mathrm{d}x\,\mathrm{d}\Omega\,\mathrm{d}t\) with isotropic \(\Omega\) and assign weights \(w_0\) so that \(\sum_p w_0=\int Q\). Particle weights will decrease exponentially along their trajectory.

\paragraph{Unlimited (reference) case.}
 Particles stream on straight lines until they scatter; along a segment of length \(\ell\), apply implicit capture \(w\leftarrow w\,\exp(-\!\int_0^\ell \sigma_a\,\mathrm{d}s)\). Candidate scatters are drawn at exponential rate \(\sigma_s\); upon a collision, resample \(\Omega\) isotropically and continue. Terminate on exit, end of step, or \(w<w_{\min}\). Track the angle-average \(\langle\Psi\rangle=(4\pi)^{-1}\int_{S^2}\Psi\,\mathrm{d}\Omega\) via a track-length estimator in space--time bins.

\paragraph{Non-scattering variant.}
Particles never change direction, but lose weight with rate  \(\sigma_t=\sigma_a+\sigma_s\), i.e., \(w\leftarrow w\,\exp(-\!\int \sigma_t\,\mathrm{d}s)\). The resulting \(\langle\Psi\rangle\) along trajectories is supplied to another solver to account for higher-order scattering.

\paragraph{\(N_s\)-scatter--limited variant.} At birth or the start of the time step each particle is assigned a scatter counter $n$. The simulation of a particle changes based on its current scatter count.

\begin{itemize}
  \item \emph{Pre-limit} (\(n<N_s\)): between collisions, use implicit capture with \(\sigma_a\) and sample next scatter at rate \(\sigma_s\); after each scatter, set \(n\leftarrow n+1\) and resample \(\Omega\).
  \item \emph{Post-limit} (\(n=N_s\)): no further collisions; stream with attenuation \(\sigma_t\) only.
\end{itemize}
Maintain two tallies: \(\langle\Psi\rangle_{\text{pre}}\) for segments with \(n\le {N_s}-1\) and \(\langle\Psi\rangle_{\text{post}}\) for the \(n={N_s}\) leg. Use \(\langle\Psi\rangle_{\text{post}}\) to form the exactly-\({N_s}\) collided source, e.g., \(Q_{{N_s}+1}=(\sigma_s/4\pi)\,\langle\Psi_{N_s}\rangle\), and iterate to recover \(({N_s}{+}1)\)-and-higher collisions. For \({N_s}=0\), this reduces to the standard uncollided/collided split.

\paragraph{Pseudocode (one time step, inline).}
Initialize \((x,\Omega,t)\sim q\), \(w\leftarrow w_0\), \(n\leftarrow 0\). While alive: if \(n<{N_s}\), sample a candidate collision at rate \(\sigma_s\), stream to the nearest of boundary/step end/collision, tally into \(\langle\Psi\rangle_{\text{pre}}\), attenuate \(w\) with \(\sigma_a\); if a collision occurs, resample \(\Omega\), set \(n\leftarrow n+1\). Else (\(n={N_s}\)): stream to termination without collisions, tally into \(\langle\Psi\rangle_{\text{post}}\), attenuate with \(\sigma_t\). Terminate on exit, step end, or \(w<w_{\min}\). Finally, pass \(\langle\Psi\rangle_{\text{post}}\) to the collided solver and iterate as needed.

\subsection{Quasi Monte Carlo (QMC) }
Quasi–Monte Carlo (QMC) methods replace pseudorandom draws on the unit cube with low-discrepancy points to reduce variance for a given sample size. In simulation, these unit-cube points are typically transformed to the required physical inputs via inverse-CDF maps (or equivalent measure-preserving transforms). In this work we use a Halton sequence because it is simple, extensible in dimension, and performs robustly for the moderate dimensions relevant here.

For transport simulation, as discussed here, a natural QMC approach is to first sample the Poisson-distributed number of scattering events to fix the stochastic dimension, and then use the low-discrepancy points to sample the 
 scattering locations/directions. However, in non-homogeneous materials, the event rate varies along the path, so 
the number of scatter events is path-dependent the required transform is not known a priori. This makes a fixed-dimension QMC mapping cumbersome and usually impractical, eroding the advantage of low-discrepancy sampling in that setting. In our hybrid approach, we see improved convergence rates by simply replacing the pseudo-random numbers of the MC algorithm with the quasi-random numbers given by a Halton sequence. We believe this to  be due to the resampling process being scatter-free, effectively limiting the required dimension of random numbers to the degrees of freedom defining such a straight-line trajectory. 

\section{Numerical examples}
\begin{figure}
    \centering
    \includegraphics[width=1\linewidth]{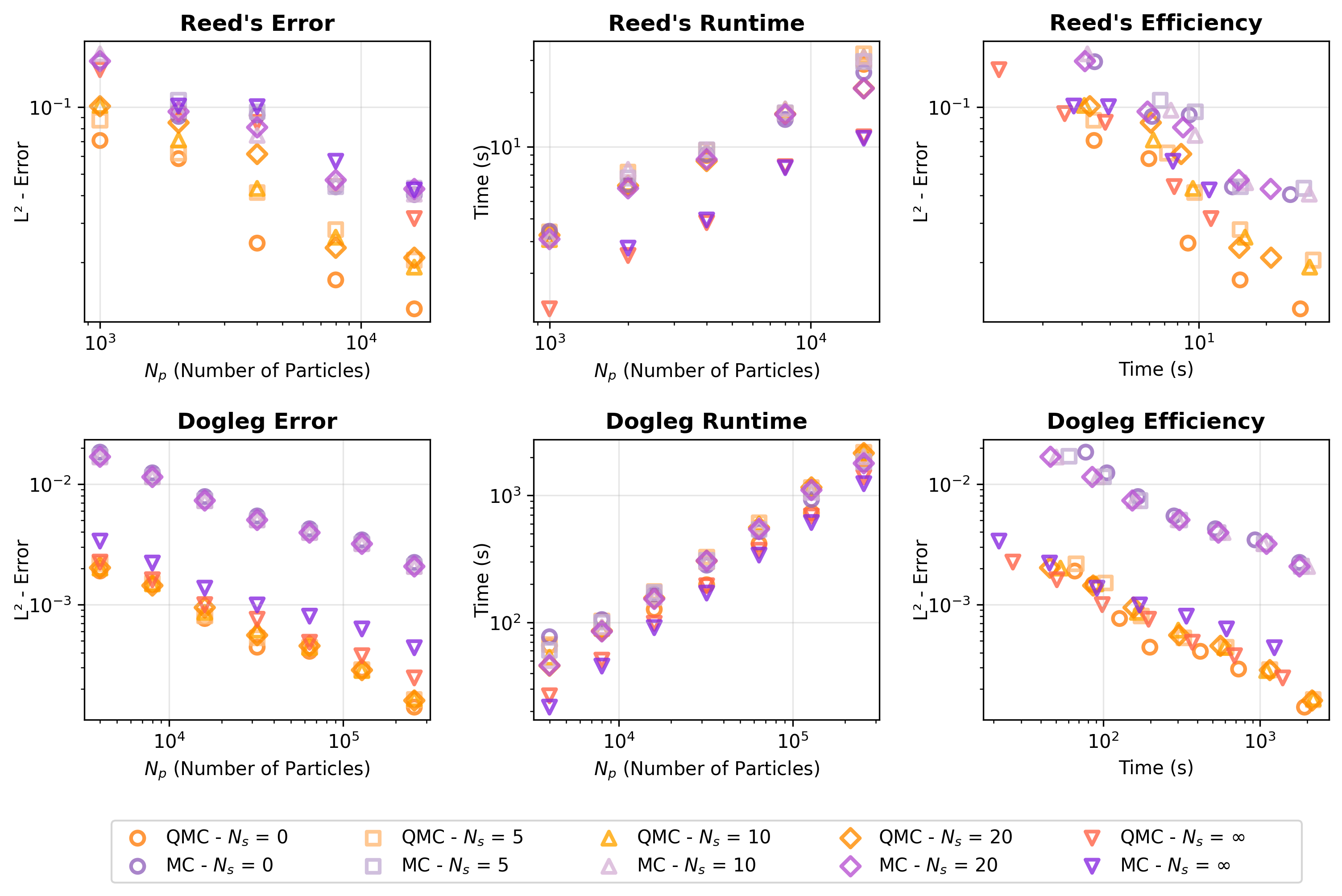}
    \caption{Summary of numerical experiments for Reed's problem (top row) and the Dogleg problem (bottom row). Depicted are the $L^2$--error as  depending on  the particle count $N_p$(left column), the runtime as depending on $N_p$ (middle column) and $L^2$--error depending on the runtime(right column) as a measure of efficiency. Runs that appear closer to the bottom left provide higher accuracy per runtime. We compare QMC runs(orange) to MC runs(purple) for various values of the scatter limit $N_s\in \{0,5,10,20\}$ as well as $N_s = \infty$(non-hybrid/regular (Q)MC run). }
    \label{fig:results}
\end{figure}

As numerical examples we include Reed's problem \cite{reed1971new}, a 1D problem with heterogenous materials and commonly used as a test case as well as a 2D problem  based on the Kobayashi Dogleg benchmark \cite{kobayashi2000nea3d}. Geometries and material properties for either problem are included in the Appendix \ref{app:geom}. Both problems are run as equilibrium problems, i.e. with a single infinite time step. While the algorithm is written for and works on explicitly time dependent problems as well, we defer to future work for numerical examples. We restrict our analysis to cases with  an $ S_4$-solver for $\Psisup{N_s}$. The solver specifics are outlined in \cite{krotz_2024}. In that work we also demonstrated that in the $N_{s} = 0$ case the order of the discrete ordinates solver does affect the overall accuracy marginally at best. Additional tests for $N_{s}>0$, not discussed here, demonstrate that this remains true here as well. All runs were performed sequentially on the same local machine using code built on top of MCDC, developed by the Center for Exascale Monte-Carlo Neutron Transport (CEMeNT) \cite{mcdc}.
Results depicted in figure \ref{fig:results} focus on the convergence in $L^2$ and run time as a function of the number of particles, number of scatter events per particle $N_{s}$ and based on whether QMC or regular MC was used. The last column in these figures depicts the $L^2$-errors as function of the runtime, to analyze overall computational efficiency: more efficient runs appear either lower (more accurate), further to the left (faster) than less efficient runs, or both. Spatial resolution was kept constant throughout with $N_x = 80$ for Reed's problem and $N_x \times N_y = 30 \times 50$ for the Dogleg benchmark. $N_{s}$ was varied between $\{0,5,10,20,\infty\}$, with $N_{s}=\infty$ corresponding to a pure MC simulation. Relabeling is unnecessary in this case and therefore skipped. 

Significant key observations in these figures are the following: The hybrid methods, independently of $N_{s}$, see significant improvements in accuracy by using QMC at no additional cost in terms of runtime. This effect is much less pronounced for the $N_{s}=\infty$, i.e. the non-hybrid case. Hybrid runs that use QMC over MC appear to be more accurate per particle than non-hybrid runs. This is solidified by improved convergence rates of QMC runs for the hybrid method, as outlined in table \ref{tab:sample}. 
This still holds for regular MC in the 1D, but not in the 2D example. It is not clear whether this is a question of dimensionality or problem dependent.  
Non-hybrid runs are, for the same number of particles, faster. This is expected as no relabeling is happening. For QMC runs the additional accuracy of hybrid runs makes up for the longer runtime, leading to QMC hybrid runs  being comparably or more efficient than non-hybrid runs. 
For hybrid runs the effect of different $N_{s}$ appears to have little effect on the run time. This is likely due to the fact that higher $N_{s}$ results in what seems to be proportional reduction in the number of iterations needed for the $S_N$-solver, also shown in table \ref{tab:sample}.

\begin{table}
    \centering
    \caption{ Comparison of number of $S_N$ iterations and convergence rate as function of allowable number of scatter events $N_s$ for Reed's problem(left) and the Dogleg problem (right).
    $S_N$ iterations are averaged over runs with various $N_p$, but $N_p$ dependence appears negligible.
    Convergence rate $\alpha$ is determined via least-square fit of $cN_p^{\alpha}$ to the $L^2$--distance between approximate and reference solutions. We distinguish $\alpha_{MC}$ and $\alpha_{QMC}$ based on the random numbers used.
   }
    \begin{tabular}{c|cll|crr}
        \toprule
                                 &                  &  Reed's problem& & & Dogleg & \\ 
        \midrule
         $N_{s}$    &  $S_N$-Iterations    &  $\alpha_{QMC}$ &  $\alpha_{MC}$&  $S_N$-Iterations &  $\alpha_{QMC}$ &  $\alpha_{MC}$ \\\midrule \midrule
                 0     & 31  & -0.69 & -0.5 &          7        &           -0.61 & -0.49 \\
         \midrule
                 5     & 26 &-0.64 &  -0.51 &          4        &           -0.61 & -0.49 \\
         \midrule
                 10    &  21.4 &  -0.63 &-0.53 &          1        &           -0.61 & -0.49 \\
         \midrule
                 20    &  12.2 &  -0.64 &-0.49 &          1        &           -0.61 & -0.49 \\
         \midrule
                 $\infty$  &  0  & -0.56  & -0.46   &          0        &           -0.53 & -.47 \\
         \midrule
        \bottomrule
    \end{tabular}
    
    \label{tab:sample}
\end{table}

\section{Conclusions}
We showed that with the addition of a relatively cheap and standard deterministic solver and a relabeling step, standard MC algorithms for particle transport can reap the benefits of that come with the use of quasi-random numbers with only minor, localized changes. Namely (i) replacing pseudorandom draws by low-discrepancy points without altering the solver’s existing inverse-CDF/transport maps, (ii) inserting a simple scatter counter with a post-limit non-scattering leg, and (iii) tallying \(\langle\Psi\rangle\) before/after the \(N_s\)-th scatter—transport solvers. This gives a plug-in pathway to QMC \emph{without} substantial changes to the current sampling procedure. This keeps code disruption low while unlocking QMC’s variance-reduction benefits inside otherwise conventional workflows.

We also introduced \emph{multi-scatter hybrids} (\(N_s>0\)), which solve a predetermined number of collisions with (Q)MC and delegate the remainder to a deterministic “clean-up’’ solve for highly scattered particles. Although the direct impact of multi-scattering on accuracy and wall time was marginal in our tests, this avenue is still compelling: We only reported runtimes in a sequential setting. However, multi-scatter QMC is far easier to parallelize than the coupled clean-up stage, promising speedups in non-sequential (e.g., many-core or distributed) runs. 

We observe, not unexpectedly, a direct trade-off between the number of iterations used in the 'clean up' and the number of allowable scatters in the (Q)MC-stage. 

A practical tuning knob emerges from this tradeoff between the number of MC-handled scatters and the iterations required by the clean-up solver: larger \(N_s\) reduces angular-coupling burden downstream, enabling cheaper deterministic settings in high-scattering regimes. This division of labor can be automated to target either fixed accuracy at minimum time or maximum accuracy under a time budget. One can envision a hybrid with a fully parallel (Q)MC leg and relatively high $N_s$, combined with a deterministic solver as cheap as a simple diffusion approach, for minimal runtimes for example. Other examples could include larger, $N_s$-dependent time steps, as with increased $N_s$ fewer particles enter the deterministic regime.

Finally, throughout this manuscript we chose to resolve the non-uniquenesss of equation \eqref{eq:uncollided} by conditioning $\Psi_n$ to solve \eqref{eq:n-split} for all $n\le N_s$ for a fixed $N_s$. This is an intuitive choice, but far from the only one, leaving tremendous space for future design choices. As a simple next step \(N_{s}\) can be made to depend on space, cross sections, time, or other indicators (e.g., optical thickness, local stiffness, or estimated variance), opening a broad path for adaptive strategies and future exploration of problem-dependent multi-scatter schedules and plug-and-play QMC integrations.

\section*{ACKNOWLEDGEMENTS}
This work was supported by the Center for Exascale Monte-Carlo Neutron Transport
(CEMeNT), a PSAAP-III project funded by the Department of Energy under Grant
DE-NA003967.

\newpage
\bibliographystyle{physor2026}
\bibliography{main}
\newpage
\appendix
\section{Uniqueness and alternative constraints} \label{app:uniqueness}
As a standalone equation, equation \eqref{eq:uncollided} can be rewritten as
\begin{align}
     \del_t (u+v)+ \Om \cdot \nabla_{\x} (u+v)+\lambda_u u + \lambda_v v= \lambdas \mathcal{S}[u] + s,  \label{eq:uv}
\end{align}
with material parameters $\lambda_u,\lambda_v$ and $\lambdas$, as well as source $s$ and appropriate boundary and initial data for $u$ and $v$. A solution to this equation is a tuple of functions $(u,v)$ that jointly fulfill the imposed conditions. Such a tuple is not determined uniquely by \eqref{eq:uv}.

For example, if we impose $v \equiv 0$ the, then unique, solution can be found by setting $u$ to be a solution to \eqref{eq:boltzmann} and $v \equiv 0$. If, on the other hand, $u\equiv 0$ was imposed $v$ solves \eqref{eq:boltzmann} with the replacements $0\rightarrow\sigmas$ and $\sigmat\rightarrow \sigmaa$.  (We are assuming compatible boundary conditions and initial data with respect to the imposed constraints.) 
All possible solutions to \eqref{eq:uv} are pointwise convex combinations of these two extreme cases. 
To enforce uniqueness an additional constraint $\bm{Law}(u,v)$ on the relation of $u,v$ is required. Examples we have already seen are \begin{align}
    \bm{Law}(u,v):&  \left\{u\equiv 0\right\} \\
    \bm{Law}(u,v):&  \left\{v\equiv 0 \right\} \\
    \bm{Law}(u,v):&  \left\{u =\sum_{n=0}^{N_s-1} \Psi_n, v=\Psi_{N_s}, \text{s.t.} \Psi_n \text{solve \eqref{eq:n-split}} \right\} 
\end{align}

Many versions of valid $\bm{Law}(u,v)$ can be found through small changes to the scatter-limited MC-algorithm. Examples include, but are not limited to, making $N_s$ spatially or material dependent.
 We posit that the freedom in choosing a specific relation $\bm{Law}(u,v)$ to force uniqueness should be seen as a feature rather than a bug.

\section{Exponential entry into the high–scatter regime}\label{app:transfer}
In the main body of the manuscript we claim that particles will transition to the high scatter regime (more than $N_s$ scatters) at an exponential rate, if $\sigmas \neq 0$, necessitating relabeling to avoid the deterministic solver determining the overall accuracy. This section intends to provide brief proof of this claim.
Assume first that the scattering rate is constant (per cell), \(\sigma_s>0\).
Then successive inter–scatter times are i.i.d.\ \(\mathrm{Exp}(\sigma_s)\).
If a particle has already scattered \(N_s-n\) times, the additional time to accrue
\(n{+}1\) scatters is
\[
T_{n+1}\sim \mathrm{Erlang}(n{+}1,\sigma_s).
\]
In particular, the transition time from source is
\(T_{N_s+1}\sim \mathrm{Erlang}(N_s{+}1,\sigma_s)\) with exact tail
\[
\mathbb{P}(T_{N_s+1}>\Delta t)
= e^{-\sigma_s \Delta t}\sum_{j=0}^{N_s}\frac{(\sigma_s \Delta t)^j}{j!}\!,
\qquad \Delta t>0.
\]
With the monotonicity $\mathbb{P}(T_{n+1}>\Delta t) > \mathbb{P}(T_{n}>\Delta t)$, this immediately yields two-sided, explicit bounds:
\[
e^{-\sigma_s \Delta t}
\;\le\; \mathbb{P}(T_{0}>\Delta t)\; \le\; \dots \;\le\;
\mathbb{P}(T_{N_s+1}>\Delta t)
\;\le\;
e^{-\mu \Delta t}\Bigl(1-\tfrac{\mu}{\sigma_s}\Bigr)^{-(N_s+1)}
\quad\text{for any } \mu\in(0,\sigma_s).
\]

Hence \(\mathbb{P}(T_{n+1}>\Delta t)=\mathcal{O}(e^{-\mu \Delta t})\) for any fixed \(\mu\in(0,\sigma_s)\):
particles enter the \(>\!N_s\) scatter regime at an \emph{exponential} rate. This constant-\(\sigma_s\) assumption
is sufficient for per-cell analysis.
If the scattering cross sections varies but is bounded
\[
0<\sigma_s^{\min}\;\le\;\sigma_s(\cdot)\;\le\;\sigma_s^{\max},
\]
then by coupling with homogeneous Poisson processes one has the stochastic sandwich
\(T_{n+1}^{(\sigma_s^{\max})}\ \preceq_{\mathrm{st}}\ T_{n+1}\ \preceq_{\mathrm{st}}\ T_{n+1}^{(\sigma_s^{\min})}\).
Consequently
\[
e^{-\sigma_s^{\max}\Delta t}
\;\le\;
\mathbb{P}(T_{n+1}>\Delta t)
\;\le\;
e^{-\mu \Delta t}\Bigl(1-\tfrac{\mu}{\sigma_s^{\min}}\Bigr)^{-(N_s+1)}
\quad\text{for any } \mu\in(0,\sigma_s^{\min}).
\]
Thus the exponential transition rate persists under \(0<\sigma_s^{\min}\le \sigma_s(\cdot)\le \sigma_s^{\max}\).

\section{Geometry   and select runs of numerical Problems}\label{app:geom}
In this section we will provide the material properties and geometries used in our numerical examples, so that a reader could reproduce these runs or use the same set ups as benchmarks for different methods. We provide outputs from some example runs as well.
\subsection{Reed's problem}
\begin{figure}[h]

        \centering
    \includegraphics[width=.8\linewidth]{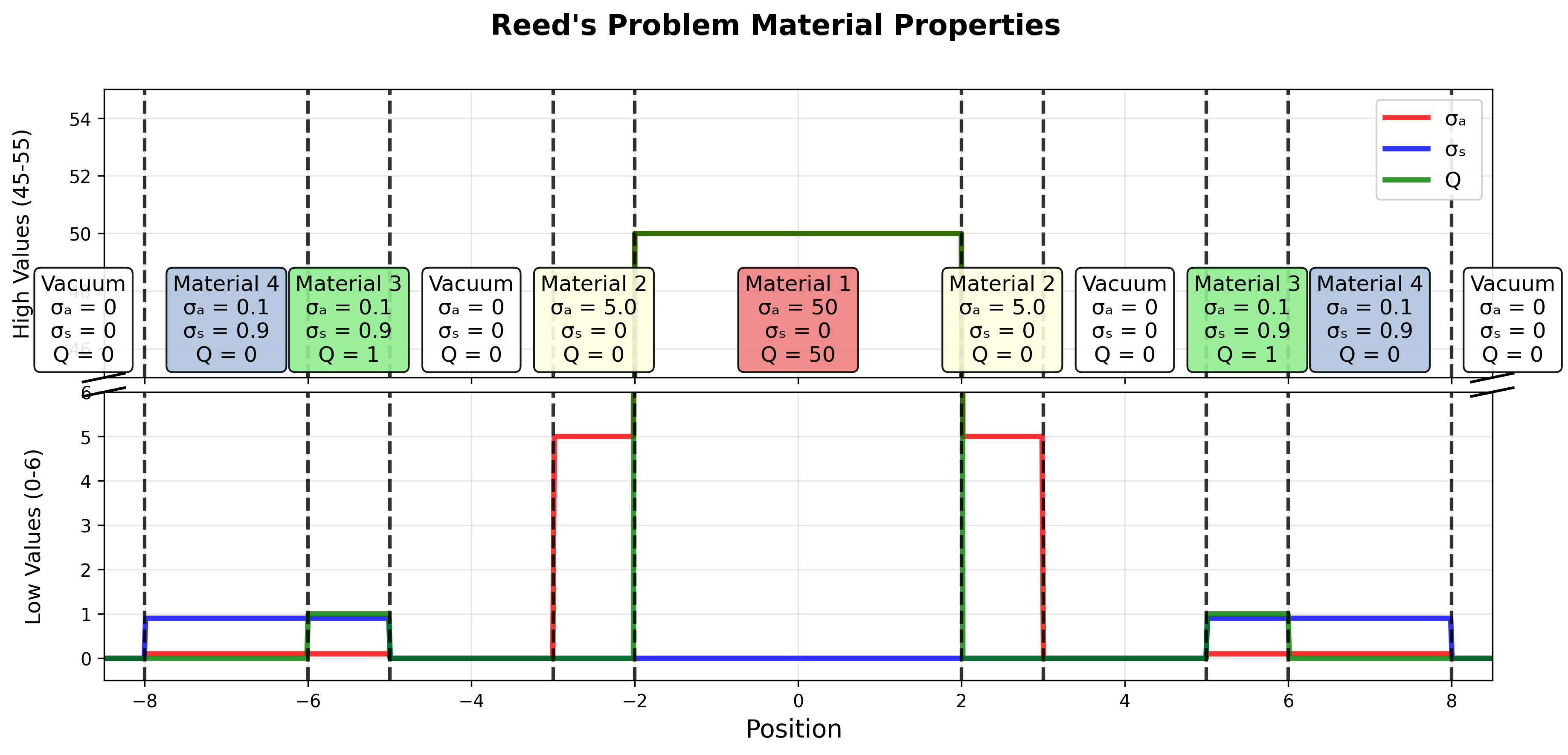}
   \caption{Material parameters and geometry for Reed's problem}
    \label{fig:reeds_mat}
\end{figure}

\begin{figure}
    \centering
    \includegraphics[width=1\linewidth]{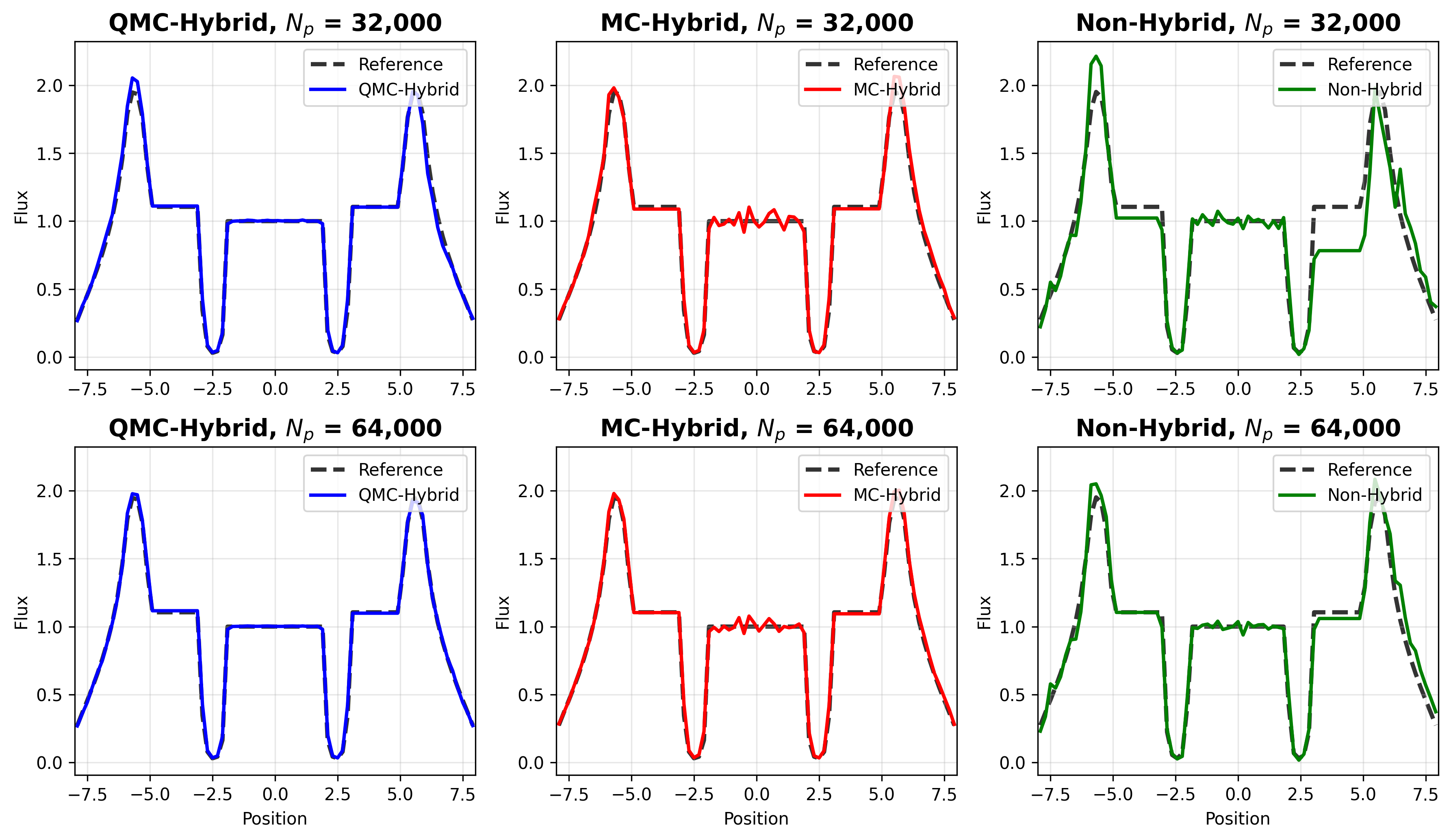}
    \caption{Sample outputs of $\langle \Psi \rangle$ for Reed's problem with $N_p = 32000$(top row) and $N_p = 64000$(bottom row). We compare results to a reference (dashed) for  QMC--hybrid runs(left), MC--hybrid runs(middle) and non-hybrid runs(right). The scattering limit $N_s$ is set to zero for hybrid runs and $\infty$ for non-hybrid runs.}
    \label{fig:reeds_runs}
\end{figure}

Reeds' problem is a 1D benchmark problem pulled from \cite{reed1971new}. It is made up of 4 slabs of different materials reflected at $x=0$. We removed the reflective boundary and replaced it via domain duplication. The layout simulated, including values for $\sigma_a,\sigma_s$ and $Q$ are depicted in figure \ref{fig:reeds_mat}. Example runs for a few different particle counts $N_p$ and different methods (QMC,MC and non-hybrid) are depicted in figure \ref{fig:reeds_runs}.
\subsection{Dogleg problem}
\begin{figure}
            \centering
            \includegraphics[width=1\linewidth]{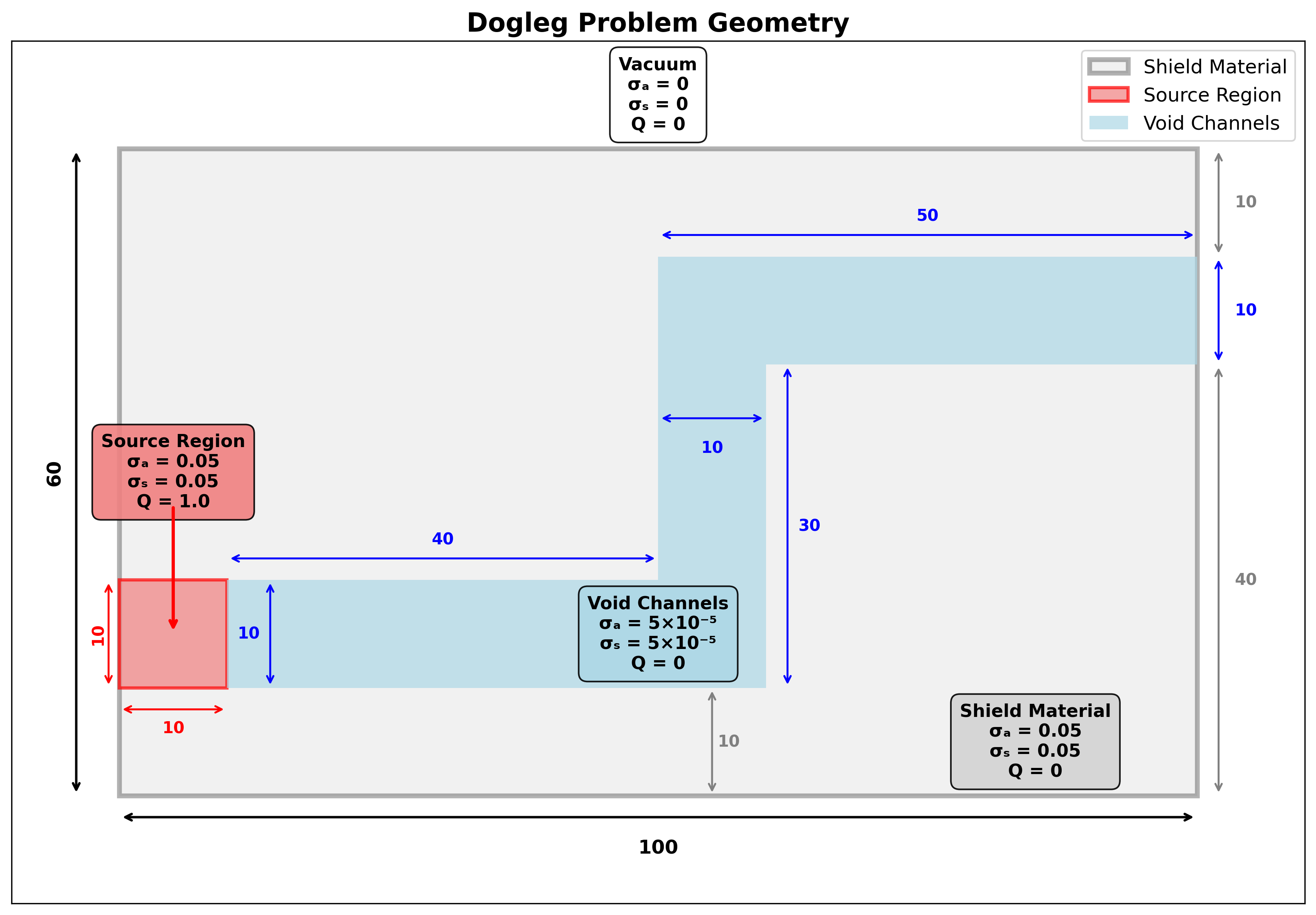}
            \caption{Material parameters and geometry for the Dogleg problem (rotated by $90^\circ$)}
            \label{fig:dogleg_mat}
        \end{figure}

\begin{figure}
            \centering
            \includegraphics[width=.9\linewidth]{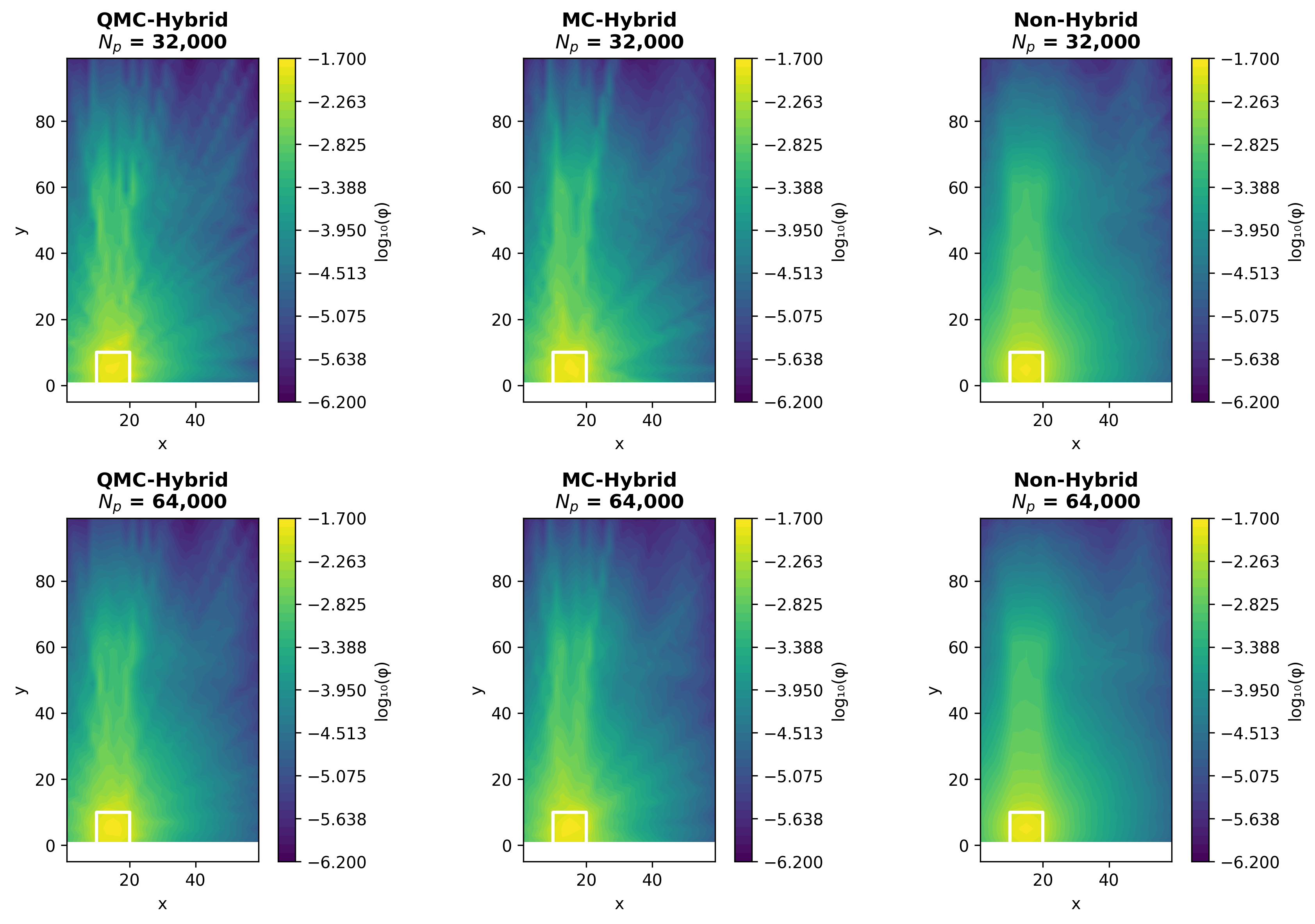}
            \caption{ Sample outputs of $\langle \Psi \rangle$ for the Dogleg problem with $N_p = 32000$(top row) and $N_p = 64000$(bottom row). We compare results for  QMC--hybrid runs(left), MC--hybrid runs(middle) and non-hybrid runs(right). The scattering limit $N_s$ is set to zero for hybrid runs and $\infty$ for non-hybrid runs.}
            \label{fig:dogleg_runs}
        \end{figure}

The Dogleg problem is a 2D problem inspired by the 3D-Dogleg problem originally suggested in \cite{kobayashi2000nea3d}. It is a simple near-vacuum duct winding through an absorbing and scattering material with a source at the start of the duct. We turned the 3D problem into a 2D problem by infinitely extending the duct and source along the z-direction. The schematic layout including values for $\sigma_a,\sigma_s$ and $Q$ throughout the domain is depicted in figure \ref{fig:dogleg_mat}. Example runs for a few different particle counts $N_p$ and different methods (QMC,MC and non-hybrid) are depicted in figure \ref{fig:dogleg_runs}.

\end{document}